\documentclass[12pt]{amsart}
\usepackage{amssymb}
\usepackage[all]{xy}
\usepackage[pdftitle={Classification of Heisenberg groups}]{hyperref}

 \rm

\newcommand{\abs}[1]{\left\lvert#1\right\rvert}

\newcommand{\C}{{\mathbb{C}}}
\newcommand{\Q}{\mathbf Q}
\newcommand{\R}{\mathbf R}
\newcommand{\Z}{\mathbf Z}

\renewcommand{\phi}{\varphi}

\newcommand{\dash}{\nobreakdash-}

\newcommand{\A}{{\mathbf{A}}}

\theoremstyle{plain}
\newtheorem*{thm*}{Theorem}

\newtheorem*{cor*}{Corollary}
\newtheorem*{thma}{Theorem~A}
\newtheorem*{thmb}{Theorem~B}
\newtheorem*{thmc}{Theorem~C}

\theoremstyle{definition}
\newtheorem*{defn*}{Definition}

\newtheorem*{rmk*}{Remark}
\newtheorem*{rmka}{Remark~A}
\newtheorem*{rmkb}{Remark~B}
\newtheorem*{rmkc}{Remark~C}
\newtheorem*{eg*}{Example}
\title[Classification of Heisenberg groups]{Decomposition of phase space and\\Classification of Heisenberg groups}

\author{Amritanshu Prasad}
\address{The Institute of Mathematical Sciences, CIT campus Taramani, Chennai 600113.}

\author{M.~K.~Vemuri}
\address{Chennai Mathematical Institute, Plot H1, SIPCOT IT Park, Padur~PO, Siruseri 603103.}


\begin{document}

\begin{abstract}
  Is every locally compact abelian group which admits a Heisenberg central extension isomorphic to the product of a locally compact abelian group and its Pontryagin dual?
  An affirmative answer is obtained for all the commonly occurring types of abelian groups having Heisenberg central extensions, including Lie groups and certain finite Cartesian products of local fields and ad\`eles.
  Furthermore, for these types of groups, it is found that the isomorphism class of the abelian group determines the Heisenberg group up to isomorphism, thereby providing a classification of such Heisenberg groups.
\end{abstract}

\maketitle


\section{Introduction}\label{sec:introduction}
In his mathematical formulation of quantum kinematics \cite[Chap.~IV, Sec.~14]{Weyl50}, Hermann Weyl introduced the concept of Heisenberg group as a special kind of central extension of an abelian group (which is thought of as phase space) by the circle group.
His formulation allowed the phase space to be quite general, and besides real vector spaces, he considered finite abelian groups.
A natural class of Heisenberg groups are the standard Heisenberg groups defined in Section~\ref{sec:defin-heis-groups}, where the phase space is the product of a locally compact abelian group (thought of as position space) with its Pontryagin dual (thought of as momentum space).
However, the theory of Heisenberg groups can be formulated independently of such a decomposition (see, for example \cite{MR1116553}).
In some important applications, the phase space associated to a Heisenberg group does not come with any natural decomposition into position and momentum spaces.
In fact, it is not known to the authors whether such a decomposition exists in general.

In this article it is shown that, in a large number of cases including many of the standard examples, such a decomposition of the phase space exists and that the Heisenberg group is isomorphic to a standard Heisenberg group.
In fact, the isomorphism class of the phase space (as a locally compact abelian group) determines the Heisenberg group up to isomorphism.
We rely on known structure theorems for certain families of locally compact abelian groups.
Since no structure theorem is available for locally compact abelian groups in general (even self-dual ones), new ideas may be needed to tackle the general case.

The organization of this paper is as follows:
In Section~\ref{sec:defin-heis-groups}, we recall the definition of Heisenberg group following \cite{MR1116553} and introduce standard Heisenberg groups.
In Section~\ref{sec:heis-lie-groups}, we treat the case where the phase space is a Lie group.
Isomorphism classes of Lie Heisenberg groups are found to be in canonical bijective correspondence with isomorphism classes of abelian Lie groups with finitely many connected components.
In Section~\ref{sec:divis-tors-free}, we treat the case where the phase space is a divisible and torsion-free locally compact abelian group.
This includes Heisenberg groups whose phase space is a finite Cartesian product of local fields of characteristic zero and ad\`eles of algebraic number fields.
Isomorphism classes of Heisenberg groups with divisible torsion-free phase space are found to be in canonical bijective correspondence with divisible torsion free groups whose identity component is a Lie group.
In Section~\ref{sec:groups-prime-expon}, we treat the case where the phase space is of prime exponent (every non-trivial element has order $p$ for a fixed prime number $p$). 
This includes Heisenberg groups whose phase space is a finite Cartesian product of local fields and ad\`eles of global function fields of characteristic $p$.
Isomorphism classes of Heisenberg groups with phase space of exponent $p$ are found to be in canonical bijective correspondence with direct products of cyclic groups of order $p$.

Our results provide a simple minimal framework within which the study of many questions involving a large class of Heisenberg groups can be carried out.
The classification of finite Heisenberg groups can be interpreted as the fact that the number of irreducible $n$\dash dimensional abelian groups of unitary rotations in ray space (see \cite[p.~276]{Weyl50}) is the same as the number of isomorphism classes of abelian groups of order $n$ (see \cite{cow}).

\section{Definition of Heisenberg Group}
\label{sec:defin-heis-groups}
Let $K$ be a locally compact abelian group, and 
$U(1) =\{ z \in \C \: :\: \abs{z} = 1 \}$.
Consider a locally compact group $G$ which is a central extension:
\begin{equation}\label{hg}
\xymatrix{
1 \ar[r] & U(1) \ar[r]^i & G \ar[r]^\pi & K \ar[r] & 0
}
\end{equation}
The commutator map $G\times G\to U(1)$ defined by $(g,h)\mapsto ghg^{-1}h^{-1}$ descends to an alternating bicharacter $e:K\times K\to U(1)$ known as the commutator form.
Thus, $e$ is a homomorphism in each variable separately, and $e(k,k)=1$ for all $k\in K$.
Let $\widehat{K}$ denote
the Pontryagin dual of $K$.  Then a homomorphism $e^\flat: K \to \widehat{K}$
is obtained by setting
$$
(e^\flat(l))(k) = e(k,l).
$$
\begin{defn*}
  If $e^\flat$ is an isomorphism, then $G$ is called a {\em Heisenberg group}.
\end{defn*}
In this case, $e$ is said to be a non-degenerate alternating bicharacter.
\begin{eg*}[Standard Heisenberg group]
  Let $A$ be a locally compact abelian group.
  Denote its Pontryagin dual by $\hat A$.
  For each $x\in A$, let $T_x$ denote the translation operator $T_xf(y)=f(y-x)$ on $L^2(A)$.
  For each $\chi \in \hat A$, let $M_\chi$ denote the modulation operator $M_\chi f(y)=\chi(y)f(y)$ on $L^2(A)$.
  Consider
  \begin{equation*}
    G=\{zT_xM_\chi\;|\;z\in U(1),\; x\in A\text{ and } \chi \in \hat A\}
  \end{equation*}
  as a subgroup of the group of unitary operators on $L^2(A)$.
  With $i(z)=zI$ for each $z\in U(1)$, $K=A\times \hat A$, and $\pi(zT_xM_\chi)=(x,\chi)$, $G$ is a Heisenberg group as in (\ref{hg}) with commutator form
  \begin{equation*}
    e((x,\chi),(x',\chi'))=\chi'(x)\chi(x')^{-1} \text{ for all } x,x'\in A \text{ and } \chi,\chi'\in \hat A.
  \end{equation*}
  We call $G$ the standard Heisenberg group associated to $A$.
\end{eg*}
In general, if $K$ is a locally compact abelian group and $c:K\times K\to U(1)$ is a $2$\dash cocycle (often called a multiplier) so that $G$ can be realized as $U(1)\times K$ with the group structure given by
\begin{equation*}
  (z,k)(w,l)=(zwc(k,l), k+l),
\end{equation*}
then $e(k,l)=c(k,l)c(l,k)^{-1}$.
It is not difficult to see that $c\mapsto e$ gives rise to an embedding of $H^2(K,U(1))$ into the group of alternating bicharacters of $K$.
Thus the data $(K,e)$, where $K$ is a locally compact abelian group and $e$ is a non-degenerate alternating bicharacter, determine the Heisenberg group $G$ up to isomorphism of central extensions.
We refer to $K$ as the phase space associated to $G$.
\begin{rmk*}
  If $x\mapsto 2x$ is an automorphism of $K$, then given an alternating bicharacter $e:K\times K\to U(1)$, we may define $c(k,l)=e(k/2,l/2)^2$ to construct a central extension with alternating bicharacter $e$.
  But if $x\mapsto 2x$ is not an automorphism of $K$, it is not clear that every non-degenerate alternating bicharacter is the commutator form of a Heisenberg group.
  When $K$ is Lie or of exponent $2$, this is a consequence of Theorems~A and C.
\end{rmk*}

\section{Lie Heisenberg groups}
\label{sec:heis-lie-groups}
A Lie group will be assumed to have countably many connected components.
\begin{thma}
  Let $K$ be an abelian Lie group and let $e:K\times K\to U(1)$ be a non-degenerate alternating bicharacter.
  Then there exists an abelian Lie group $A$ and an isomorphism $\phi:A\times \hat A\to K$ such that
  \begin{equation*}
    e(\phi(x,\chi),\phi(x',\chi'))=\chi'(x)\chi(x')^{-1} \text{ for all } x,x'\in A\text{ and }\chi,\chi'\in \hat A.
  \end{equation*}
\end{thma}
\begin{rmka}
  In the above theorem, the isomorphism class of $A$ is not uniquely determined by the isomorphism class of $(K,e)$.
  For example, $A=\Z$ and $A=U(1)$ give rise to isomorphic Heisenberg groups.
  However, $A$ is uniquely determined up to isomorphism if we stipulate that $A$ have only finitely many connected components.
\end{rmka}
\begin{proof}
  $K$ can be decomposed as
  \begin{equation*}
    K=K_f\times K_\R\times K_t \times K_\Z,
  \end{equation*}
  where $K_f$ is a finite abelian subgroup, $K_\R$ is a subgroup isomorphic to the additive group of a finite dimensional real vector space, $K_t$ is a subgroup isomorphic to a compact real torus, and $K_\Z$ is a finitely generated free abelian subgroup.
  Thus $\hat K$ inherits a decomposition
  \begin{equation*}
    \hat K = \hat K_f \times \hat K_\R\times \hat K_t\times \hat K_\Z,
  \end{equation*}
  where $\hat K_f$, $\hat K_\R$, $\hat K_t$ and $\hat K_\Z$ are subgroups of $\hat K$ isomorphic to the Pontryagin duals of $K_f$, $K_\R$, $K_\Z$ and $K_t$ respectively.
  The isomorphism $e^\flat$ can be written as a matrix $(e_{ij})$, where $e_{ij}:K_j\to \hat K_i$ for $i,j\in \{f, \R, t, \Z\}$.
  Many of these entries are forced to be trivial, as there are no non-trivial homomorphisms from some types of groups to others.
  We see immediately that $e^\flat$ has the form:
  \begin{equation*}
    e^\flat=
    \begin{pmatrix}
      e_{ff} & 0 & 0 & e_{f\Z}\\
      0 & e_{\R\R} & 0 & e_{\R\Z}\\
      e_{tf} & e_{t\R} & e_{tt} & e_{t\Z}\\
      0 & 0 & 0 & e_{\Z\Z}
    \end{pmatrix}
    .
  \end{equation*}
  Noting that $(e^\flat)^{-1}$ must also have a similar matrix form, we conclude that $e_{ii}$ must be an isomorphism for each $i\in \{f,\R,t,\Z\}$.
  Given an automorphism $\alpha$ of $K$, let $e^\alpha (k,l)=e(\alpha(k),\alpha(l))$.
  The automorphism $\alpha$ itself can be represented by a matrix $(\alpha_{ij})$ of a similar form, with $\alpha_{ii}$ an automorphism of $K_i$.
  We have
  \begin{equation}
    \label{eq:4}
    (e^\alpha)^\flat = \alpha^* e^\flat \alpha.
  \end{equation}
  Here $\alpha^*: \hat K\to \hat K$ is the transpose of $\alpha$, defined by $\alpha^*(\xi)(k)=\xi(\alpha(k))$ for $k\in K$, $\xi\in \hat K$.
  For the purpose of proving Theorem~A, it is always permissible to replace $e$ by $e^\alpha$.
  
  Since $e_{ii}$ is an isomorphism for each $i\in \{f,\R,t\}$, we may take $\alpha_{i\Z}:K_\Z\to K_i$ such that $e_{i\Z}=e_{ii}\alpha_{i\Z}$ for each $i\in \{f,\R\}$,
  and $\alpha_{t\Z}:K_\Z\to K_t$ such that
  \begin{equation*}
    e_{tt}\alpha_{t\Z}=-\alpha^*_{f\Z}e_{f\Z}-\alpha^*_{\R\Z}e_{\R\Z}+e_{t\Z}-\alpha^*_{t\Z}e_{\Z\Z}.
  \end{equation*}
  Consider
  \begin{equation}
    \label{eq:3}
    \alpha=
    \begin{pmatrix}
      1 & 0 & 0 & -\alpha_{f\Z}\\
      0 & 1 & 0 & -\alpha_{\R\Z}\\
      0 & 0 & 1 & -\alpha_{t\Z}\\
      0 & 0 & 0 & 1
    \end{pmatrix}
    .
  \end{equation}
  Its transpose, in matrix form, is
  \begin{equation}
    \label{eq:2}
    \alpha^*=
    \begin{pmatrix}
      1 & 0 & 0 & 0\\
      0 & 1 & 0 & 0\\
      -\alpha^*_{f\Z} & -\alpha^*_{\R\Z} & 1 & -\alpha^*_{t\Z}\\
      0 & 0 & 0 & 1
    \end{pmatrix}
  \end{equation}
  Furthermore, the fact that $e$ is alternating means that $(e^\flat)^*(l)(k)=e^\flat(k)(l)=e(l,k)=e(k,l)^{-1}=e(-k,l)=e^\flat(l)(-k)$.
  Consequently, 
  \begin{equation}
    \label{eq:1}
    \begin{pmatrix}
      e^*_{ff} & 0 & 0 & e^*_{tf}\\
      0 & e^*_{\R\R} & 0 & e^*_{t\R}\\
      e^*_{f\Z} & e^*_{\R\Z} & e^*_{\Z\Z} & e^*_{t\Z}\\
      0 & 0 & 0 & e^*_{tt}
    \end{pmatrix}
    =
    \begin{pmatrix}
      -e_{ff} & 0 & 0 & -e_{f\Z}\\
      0 & -e_{\R\R} & 0 & -e_{\R\Z}\\
      -e_{tf} & -e_{t\R} & -e_{tt} & -e_{t\Z}\\
      0 & 0 & 0 & -e_{\Z\Z}
    \end{pmatrix}
  \end{equation}
  Using (\ref{eq:4}), (\ref{eq:2}) and (\ref{eq:1}), and by our choice of $\alpha$ in (\ref{eq:3}), we see that
  \begin{equation}
    \label{eq:5}
    (e^\alpha)^\flat=
    \begin{pmatrix}
      e_{ff} & 0 & 0 & 0\\
      0 & e_{\R\R} & 0 & 0\\
      0 & 0 & e_{tt} & 0\\
      0 & 0 & 0 & e_{\Z\Z}
    \end{pmatrix}
    .
  \end{equation}
  Therefore, it suffices to prove Theorem~A in the cases where $K=K_f$, $K=K_\R$ and $K=K_t\times K_\Z$ separately.
  \subsection*{Analysis of the finite part}
  The argument from \cite{iafhg} is repeated here for the sake of completeness.
  Choose a basis $x_1,\ldots,x_n$ of $K_f$ such that the order of $x_i$ is $d_i$ with $d_n|d_{n-1}|\cdots|d_1$.
  Let $\zeta$ denote a fixed primitive $d_1$\dash th root of unity.
  For each $i,j\in \{1,\ldots,n\}$, there exists a non-negative integer $q_{ij}$ such that
  \begin{equation*}
    e(x_i,x_j)=\zeta^{q_{ij}}.
  \end{equation*}
  Observe that $q_{ij}$ is determined modulo $d_1$ and is constrained to be divisible by $d_1/(d_i,d_j)$.
  Let $Q(e)$ denote the matrix $(q_{ij})$ associated to $e$ as above.
  
  Any automorphism $\alpha_f$ of $K_f$ is of the form
  \begin{equation*}
    \textstyle \alpha_f\big(\sum_{j=1}^n a_jx_j\big) = \sum_{i=1}^n \sum_{j=1}^n \alpha_{ij} a_j x_i
  \end{equation*}
  for a suitable matrix $(\alpha_{ij})$ of integers,
  where $\alpha_{ij}$ is determined modulo $d_i$ and is constrained to be divisible by $d_i/(d_i,d_j)$.
  Furthermore,
  \begin{equation*}
    Q(e^\alpha)=\alpha^* Q(e) \alpha.
  \end{equation*}
  In the above expression, $\alpha$ is identified with its matrix $(\alpha_{ij})$, so that $\alpha^*$ is identified with the transpose matrix.
  In particular, $Q(e^\alpha)$ can be obtained from $Q(e)$ by a sequence of column operations followed by the corresponding sequence of row operations.
  The following operations are permitted:
  \begin{enumerate}
  \item [$\beta_{i\sigma}$: ] multiplication of the $i$th row and $i$th column by an integer $\sigma$ such that $(\sigma,d_i)=1$.
  \item [$\pi_{ij}$: ] interchange of $i$th and $j$th rows and columns when $d_i=d_j$.
  \item [$\alpha_{ij\sigma}$: ] addition of $\sigma$ times $i$th row to the $j$th row, and $\sigma$ times the $i$th column to the $j$th column when $d_i/(d_i,d_j)$ divides $\sigma$.
  \end{enumerate}
  
  Since $e$ is alternating, $Q(e)$ can be taken to be skew-symmetric with zeroes along the diagonal.
  Since $e$ is non-degenerate, there must be an integer in the first row of $Q(e)$ that is coprime to $d_1$.
  First, $\beta_{i\sigma}$ for appropriate $\sigma$ can be used to make this entry equal to $1$.
  If this entry is in the $j$th column, then the constraints on $q_{ij}$ force $d_j=d_1$.
  Now, $\pi_{2j}$ transforms the top-left corner of the matrix into $\left(\begin{smallmatrix} 0 & 1 \\ -1 & 0 \end{smallmatrix}\right)$.
  In particular $d_1=d_2$.
  Finally, $\alpha_{1j\sigma}$ and $\alpha_{2j\sigma}$ can used to make the remaining entries of the first two rows and columns vanish.
  
  The form $e$ identifies the Pontryagin dual of the subgroup generated by the first (transformed) generator with the subgroup generated by the second one.
  After splitting off the subgroup spanned by these two generators, we are left with a finite abelian group with $n-2$ generators.
  By induction on $n$, we conclude that there exists a finite abelian group $A_f$, together with an isomorphism $\phi_f: A\times \hat A \to K_f$ such that
  \begin{equation*}
    e(\phi_f(x,\chi),\phi_f(x'\chi'))=\chi'(x)\chi(x')^{-1} \text{ for all } x,x'\in A\:\chi,\chi'\in \hat A.
  \end{equation*}
  \subsection*{Analysis of the real part} Since $K_\R$ is the additive group of a finite dimensional vector space, $e(k,l)=e^{iQ(k,l)}$ for all $k,l\in K_\R$ with $Q:K_\R\times K_\R\to \R$ a symplectic form.
  It is well known in this case that $K_\R$ must be an even dimensional vector space, and that $Q$ can be transformed into the standard symplectic form by an automorphism of $K_\R$.
  Hence there exists a vector space $A_\R$ and an isomorphism $\phi_\R:A_\R\times \hat A_\R\to K_\R$ such that
  \begin{equation*}
    e(\phi(x,\chi),\phi(x',\chi'))=\chi'(x)\chi(x')^{-1} \text{ for all } x,x'\in A_\R\; \chi,\chi'\in \hat A_\R.
  \end{equation*}
  \subsection*{Analysis of the toric and free parts}
  By (\ref{eq:5}), we may assume that $e_{t\Z}=0$.
  Note that $e_{\Z\Z}$ is an isomorphism of $K_\Z$ onto the Pontryagin dual of $K_t$.
  When $A_t=K_t$, and $\phi_t:A_t\times \hat A_t\to K_t\times K_\Z$ is $\mathrm{id}\times e_{\Z\Z}^{-1}$,
  we have
  \begin{equation*}
    e(\phi_t(x,\chi),\phi_t(x',\chi'))=\chi'(x)\chi(x')^{-1} \text{ for all } x,x'\in A,\:\chi,\chi'\in \hat A.
  \end{equation*}
\end{proof}
\section{Divisible torsion-free phase space}
\label{sec:divis-tors-free}
\begin{thmb}
  Let $K$ be locally compact, divisible and torsion-free, and let $e:K\times K\to U(1)$ be a non-degenerate alternating bicharacter.
  Then there exists a locally compact divisible torsion-free abelian group $A$ and an isomorphism $\phi:A\times \hat A\to K$ such that
  \begin{equation*}
    e(\phi(x,\chi),\phi(x',\chi'))=\chi'(x)\chi(x')^{-1} \text{ for all } x,x'\in A\text{ and }\chi,\chi'\in \hat A.
  \end{equation*}
\end{thmb}
\begin{rmkb}
  As in the case of Lie Heisenberg groups, $K$ does not determine $A$ up to isomorphism in Theorem~B. For example, $A=\Q$ and $A=\hat \Q$ give rise to isomorphic Heisenberg groups.
However, $A$ is uniquely determined if we stipulate that every connected component of $A$ is a Lie group.
\end{rmkb}
\begin{proof}
  According to \cite[25.33]{MR551496}, $K$ can be decomposed as 
  \begin{equation*}
    K=K_\R \times K_\A \times K_\Q\times K_s,
  \end{equation*}
  where $K_\R$ is a subgroup isomorphic to the additive group of a finite dimensional real vector space,  $K_\A$ is a restricted direct product (over all primes $p$) of finite dimensional vector spaces over $p$\dash adic fields, $K_\Q$ is a subgroup isomorphic to the additive group of a finite dimensional vector space over the rational numbers with the discrete topology, and $K_s$ is a subgroup isomorphic to a finite Cartesian power of $\hat\Q$.
  The finite dimensionality of the vector spaces over $\Q_p$ is not part of the classification of divisible torsion-free groups in \cite{MR551496}, where the minimal divisible extension of an arbitrary Cartesian power of $\Z_p$ is allowed, but is a consequence of the fact that $K$ is isomorphic to its own Pontryagin dual.
  
  Likewise, $\hat K$ can be decomposed as
  \begin{equation*}
    \hat K = \hat K_\R\times \hat K_\A\times \hat K_\Q \times \hat K_s,
  \end{equation*}
  where $\hat K_\R$, $\hat K_\A$, $\hat K_\Q$ and $\hat K_s$ are subgroups of $\hat K$ isomorphic to the Pontryagin duals of $K_\R$, $K_\A$, $K_s$ and $K_\Q$ respectively.
  As in Section~\ref{sec:heis-lie-groups}, the isomorphism $e^\flat$ can be written as a matrix some of whose entries are forced to be $0$.
  It is not difficult to see that $e^\flat$ must have the form:
  \begin{equation*}
    e^\flat =
    \begin{pmatrix}
      e_{\R\R} & 0 & 0 & e_{\R s}\\
      0 & e_{\A\A} & 0 & e_{\A s}\\
      e_{\Q\R} & e_{\Q\A} & e_{\Q\Q} & e_{\Q s}\\
      0 & 0 & 0 & e_{ss}
    \end{pmatrix}
    ,
  \end{equation*}
  where $e_{ij}$ is a homomorphism $K_j\to \hat K_i$.
  Proceeding as in Section~\ref{sec:heis-lie-groups}, we reduce Theorem~B to the cases where $K=K_\R$, $K=K_\A$ and $K=K_\Q\times K_s$.
  The first case, the analysis of the real part, has already been dealt with in Section~\ref{sec:heis-lie-groups}, and the third case is similar to the analysis of the toric and free parts treated there.
  It remains to address the case where $K=K_\A$.
  We know that $K_\A$ is a restricted direct product:
  \begin{equation*}
    K_\A = {\prod_p}' K_p,
  \end{equation*}
  where, for each prime $p$, $K_p$ is a finite dimensional vector space over $\Q_p$.
  For $p\neq q$, there are no non-trivial homomorphisms $K_p\to \hat K_q$.
  Therefore $e_{\A\A}$, when restricted to $K_p$, gives rise to an isomorphism $e_p:K_p\to \hat K_p$, and we are reduced the case where $K=K_p$ is a finite dimensional vector space over $\Q_p$.
  Let $\psi$ be a non-trivial character $\psi\in \hat \Q_p$.
  Then we may write $e(k,l)=\psi(Q(k,l))$ for all $k,l\in K_p$, where $Q:K_p\times K_p\to \Q_p$ is a non-degenerate alternating bilinear form.
  As in the real case, this form can be transformed into the standard symplectic from by an automorphism of $K_p$, and the result follows.
\end{proof}
\section{Phase space of prime exponent}
\label{sec:groups-prime-expon}
\begin{thmc}
  Let $p$ be a prime number and suppose that $K$ is a locally compact abelian group of exponent $p$.
  Let $e:K\times K\to U(1)$ be a non-degenerate alternating bicharacter.
  Then there exists a locally compact abelian group $A$ of exponent $p$ and an isomorphism $\phi:A\times \hat A\to K$ such that 
  \begin{equation*}
    e(\phi(x,\chi),(x',\chi'))=\chi'(x)\chi(x')^{-1} \text{ for all } x,x'\in A\text{ and } \chi,\chi'\in \hat A.
  \end{equation*}
\end{thmc}
\begin{rmkc}
  Once again, $K$ does not determine $A$ up to isomorphism, unless we stipulate that $A$ is a Cartesian power of a cyclic group of order $p$ with the product topology.
  This will be evident from the proof of Theorem~C.
\end{rmkc}
\begin{proof}
  Using a result of Braconnier (see \cite[25.29]{MR551496}) and the fact that $K$ is self-dual, we see that there exists a unique cardinal number $m$ such that $K$ can be decomposed as 
  \begin{equation*}
    K=K_c\times K_d,
  \end{equation*}
  where $K_c$ is a direct product of $m$ copies of $\Z/p\Z$, a compact group under the product topology, and $K_d$ is the direct sum (in the category of abelian groups) of $m$ copies of $\Z/p\Z$, and is endowed with the discrete topology.
  Therefore, $\hat K$ has a decomposition
  \begin{equation*}
    \hat K = \hat K_c \times \hat K_d, 
  \end{equation*}
  where $\hat K_c$ and $\hat K_d$ are subgroups of $\hat K$ isomorphic to the Pontryagin duals of $K_d$ and $K_c$ respectively.
  Write $e^\flat$ as a matrix
  \begin{equation}
    \label{eq:7}
    e^\flat = 
    \begin{pmatrix}
      e_{cc} & e_{cd}\\
      e_{dc} & e_{dd}
    \end{pmatrix}
    ,
  \end{equation}
  where $e_{ij}$ is a homomorphism $K_j\to \hat K_i$.
  The image of $e_{dc}$, being the continuous image of a compact group in a discrete one, is finite and may be thought of as a finite dimensional vector space over $\Z/p\Z$.
  Any lift of a basis extends to a homomorphism $s:\mathrm{Image}(e_{dc})\to K_j$.
  Let $F=s\circ e_{dc}(K_j)$.
  Then $K_c=\ker e_{dc}\times F$.
  Now let $K_c'=\ker e_{dc}$ and let $K_d'=F\times K_d$.
  Suppose that, with respect to the decompositions $K=K_c'\times K_d'$ and $\hat K=\hat K_c'\times \hat K_d'$ (where $\hat K_c'$ and $\hat K_d'$ are isomorphic to the Pontryagin duals of $K_d'$ and $K_c'$ respectively), $e^\flat$ has the matrix
  \begin{equation*}
    e^\flat = 
    \begin{pmatrix}
      e_{cc}' & e_{cd}'\\
      e_{dc}' & e_{dd}'
    \end{pmatrix}
    ,
  \end{equation*}
  then $e_{dc}'=0$.
  Now Theorem~C can be obtained as in the analysis of the real and toric parts in Section~\ref{sec:heis-lie-groups}.
\end{proof}
\subsection*{Acknowledgments}
The authors thank P.~P.~Divakaran and George Willis for their help and encouragement.

\bibliographystyle{amsplain}
\bibliography{refs}

\end{document}